\documentclass{article}
\usepackage{emsnewsletter}

\usepackage[utf8]{inputenc} 
\usepackage[T1]{fontenc}    
\usepackage{hyperref}       
\usepackage{url}            
\usepackage{booktabs}       
\usepackage{amsfonts}       
\usepackage{nicefrac}       
\usepackage{microtype}      
\usepackage{xcolor}         
\usepackage{amsmath}
\usepackage{amssymb}
\usepackage{mathtools}
\usepackage{amsthm}
\usepackage{makecell}
\usepackage[font=small]{caption}
\usepackage{multirow}
\usepackage{makecell}
\usepackage{algorithm}
\usepackage{algorithmic}
\usepackage{lscape}
\usepackage{wrapfig}

\allowdisplaybreaks
\usepackage{colortbl}
\definecolor{bgcolor}{rgb}{0.8,1,1}
\definecolor{bgcolor2}{rgb}{0.8,1,0.8}
\usepackage{threeparttable}

\def\bk{\bar \kappa}
\usepackage[capitalize,noabbrev]{cleveref}

\theoremstyle{plain}
\newtheorem{assumption}{Assumption}

\usepackage[]{todonotes}

\def\R{\mathbb{R}}

\def\R{\mathbb R}

\def\e{\varepsilon}

\newcommand{\argmin}{\operatornamewithlimits{argmin}}

\newcommand{\ls}{\left(}
\newcommand{\rs}{\right)}

\newcommand{\lb}{\left\lbrace}
\newcommand{\rb}{\right\rbrace}
\newcommand{\la}{\left\langle}
\newcommand{\ra}{\right\rangle}

\newcommand{\eqdef}{\vcentcolon=}

\providecommand{\keywords}[1]
{
  \small	
  \textbf{\textit{Keywords---}} #1
}

\def\<#1,#2>{\langle #1,#2\rangle}

\author{Dmitry Kamzolov (Mohamed bin Zayed University of Artificial Intelligence, Abu Dhabi, UAE), 
Alexander Gasnikov (MIPT, Moscow, Russia; IITP RAS, Moscow, Russia; Caucasus Mathematical Center, Adyghe State University, Maikop, Russia),
Pavel Dvurechensky (WIAS, Berlin, Germany),
Artem Agafonov (Mohamed bin Zayed University of Artificial Intelligence, Abu Dhabi, UAE; MIPT, Moscow, Russia), 
Martin Tak\'a\v{c} (Mohamed bin Zayed University of Artificial Intelligence, Abu Dhabi, UAE)}
\title{Exploiting higher-order derivatives in convex optimization methods}

\begin{document}
   
\maketitle

\section{Introduction}\label{intro}
It is well known since the works of I.~Newton \cite{newton1967methodus} and L.~Kantorovich \cite{kantorovich1949newton} that the second-order derivative of the objective function can be used in numerical algorithms for solving optimization problems and nonlinear equations and that such algorithms have better convergence guarantees. Higher-order derivatives can be efficiently used for solving  nonlinear equations as was shown by P.~Chebyshev \cite{chebyshev2018full,evtushenko2013pth}. For optimization problems, the basic idea known at least since 1970’s \cite{hoffmann1978higher} and developed in later works \cite{schnabel1984tensor,schnabel1991tensor,bouaricha1997tensor} is to approximate, at each iteration of an algorithm,  the objective by its Taylor polynomial at the current iterate, optionally add a regularization, and minimize this Taylor approximation to obtain the next iterate. In this way, the first Taylor polynomial leads to first-order methods that are very well understood, see, e.g., \cite{nemirovskij1983problem,nesterov2018lectures}, with the optimal methods existing since 1980's \cite{nemirovskij1983problem,nesterov1983method}.The second-order methods are using the second Taylor polynomial. The most famous representative is the \textit{Newton's method} that minimizes at each iteration the second-order quadratic approximation of the objective function. If the second derivative is Lipschitz continuous, the objective is strongly convex, and the starting point is sufficiently close to the solution, then this algorithm has very fast quadratic convergence and requires $\log \log \varepsilon^{-1}$ iterations to reach an $\varepsilon$-solution in terms of the objective value \cite{kantorovich1949newton}. This bound is optimal \cite{nemirovskij1983problem} even for univariate optimization problems with the possibility of using in algorithms derivatives of any order. Different modifications of the basic algorithm, such as the Damped Newton's method or the Levenberg--Marquardt algorithm achieve global convergence, but have a slower, i.e., linear, convergence \cite{polyak1987introduction,nocedal1999numerical}.
Second-order methods played also the central role in the development by Yu.~Nesterov and A.~Nemirovskii of interior point methods that have global linear convergence rate and allow proving polinomial solvability of a large class of convex optimization problems \cite{nesterov1994interior}. This theory was based on the analysis of the Damped Newton's method for the class of self-concordant functions that, in particular, include functions without Lipschitz derivatives. 

An important idea that eventually led to the current developments of tensor methods was the cubic regularization of the Newton's method, which dates back to the work of A.~Griewank \cite{griewank1981modification}. The global performance of the Cubic regularized Newton's method was analysed in 2006 by Yu.~Nesterov and B.~Polyak in \cite{nesterov2006cubic} for a large list of settings with the main assumption that the second  derivative is Lipschitz continuous. The main idea is based on the fact that, for functions with Lipshitz second derivative, the objective's model consisting of the second Taylor polynomial regularized by the cube of the norm with sufficiently large regularization parameter is an upper bound for the objective function. This model is minimized in each iteration of the algorithm, which, in particular, allowed to obtain global convergence rate $1/k^2$ for minimizing convex functions with Lipschitz second derivative. Moreover, the authors showed that the complexity of minimizing the third-order polynomial model of the objective is of the same order as the complexity of the standard Newton's method step. The Cubic regularized Newton's method was further accelerated in \cite{nesterov2008accelerating} to reach $1/k^3$ convergence rate for convex problems with Lipschitz second derivative. In 2012 R.~Monteiro and B.~Svaiter  proposed in  \cite{monteiro2013accelerated} a very perspective \textit{accelerated-proximal envelope} (see also the work \cite{carmon2022recapp}) that allowed them to develop even faster second-order method with the convergence rate $1/k^{3.5}$ for minimizing convex objectives with Lipschitz second derivative.

Another important step in the development of tensor methods was made by M.~Baes in 2009 \cite{baes2009estimate}, where the Newton--Nesterov--Polyak algorithm was generalized to the setting of convex minimization with the objective having Lipschitz $p$-th derivative for $p\ge 2$. 
The idea  was to construct a $(p+1)$-th order polynomial model that upper bounds the objective function by taking the $p$-th Taylor polynomial of the objective and regularizing it by the $(p+1)$-th power of the Euclidean norm with sufficiently large regularization parameter. The author showed $1/k^p$ global convergence rate for methods that minimize such model in each iteration and proposed an accelerated version with the rate $1/k^{p+1}$, all under the assumption of Lipschitz $p$-th derivative.
As in the world of first-order methods, where the interest in optimal methods was one of the central driving forces in 2000--2020, a natural question arose on what are the lower bounds for second- and higher-order methods and which algorithms are optimal. The results in \cite{arjevani2019oracle,agarwal2018lower,gasnikov2018hypothesis} gave the answer that the lower bound is $1/k^{(3p+1)/2}$, revealing the gap when $p\geq 3$ between the rates of existing methods an lower bounds. One of the drawbacks of tensor methods at this stage was that each iteration required to minimize a higher-order polynomial that may not be convex, leading to the high cost of each iteration and impracticality of such methods.


In 2018, a breakthrough was made by Yu. Nesterov~\cite{nesterov2021implementable} in understanding the significance of higher-order methods in modern convex optimization when the $p$-th derivative satisfies the Lipschitz condition. The breakthrough involved increasing the regularization parameter in M. Baes's approach for the regularized $p$-th Taylor polynomial of the objective. It was demonstrated that this modified Taylor model is convex. Additionally, it was shown that the convex Taylor model can be efficiently minimized for $p=3$, with similar complexity to the standard Newton's method step. Furthermore, it was proven that minimizing the convex Taylor model doesn't require the calculation of the entire tensor of third derivatives. Instead, directional third derivatives can be calculated, for example, through automatic differentiation. Lastly, a lower bound of $1/k^{(3p+1)/2}$ was established for convex problems with Lipschitz $p$-th derivative. For each $p\geq 2$, worst-case functions within the class of convex functions with Lipschitz $p$-th derivatives were constructed, and a research program was outlined to develop optimal tensor methods~\cite{nesterov2018lectures}.
It followed from the obtained lower bounds for convex optimization problems that the Monteiro--Svaiter second-order method is optimal up to a logarithmic factor since it used a complicated line-search procedure. Based on the Monteiro--Svaiter method and the convex Taylor model of the objective, in \cite{gasnikov2018substantiation,gasnikov2017universal,gasnikov2019optimal} and independently in \cite{bubeck2019near,jiang2019optimal}, the authors proposed near-optimal tensor methods for convex problems with Lipschitz $p$-th derivatives, which have the rate $1/k^{(3p+1)/2}$ up to logarithmic factors.  
In Section~\ref{optimal}, the Monteiro--Svaiter approach with Nesterov's implementable tensor steps is described.
As it was mentioned in \cite{nesterov2018lectures,doikov2020inexact}, the auxiliary line-search procedures that lead to logarithmic factors in the convergence rates of the near-optimal tensor methods may also significantly slow down the convergence in practice.
In 2022, two independent works \cite{kovalev2022first_} and \cite{carmon2022optimal} proposed optimal tensor methods with convergence rates without additional logarithmic factors. At the same time, such methods were proposed also for monotone variational inequalities \cite{adil2022line,lin2022perseus} under higher-order regularity of the operator.

The developments in the theory of tensor methods had also an important impact on the theory of second-order methods. In 2020, Yu.~Nesterov proposed an implementation of the third-order method using inexact third derivative constructed via finite differences of gradients \cite{nesterov2021superfast}. In other words, it appeared to be possible to implement the third-order method using only first and second derivatives, which lead to <<superfast>> second-order method with the convergence rate $1/k^{4}$ violating the lower bound thanks to additional assumption that the third derivative is Lipschitz continuous.
The key observation was that in the class of functions with Lipschitz second derivative the worst-case function \cite{nesterov2021implementable} for second-order methods does not have Lipschitz third derivative.
These results were further improved in \cite{kamzolov2020near,nesterov2021inexact}, where second-order algorithms with the rate $1/k^{5}$ corresponding to optimal third-order methods were proposed for minimizing convex functions with Lipschitz third derivative.
These developments are in large contrast to first-order methods, for which the worst-case function is quadratic \cite{nesterov2003introductory} and has Lipschitz derivatives of any order, thus, preventing improvements under additional assumptions. 
This line of research was continued in \cite{ahookhosh2021high}, where such reductions were shown to be possible for higher-order methods. The ideas used in superfast methods are described in Section~\ref{hyperfast}.

Tensor methods remain a very active area of research with many extensions of the existing methods. In particular, there are adaptive variants for the setting when the Lipschitz constant is not known \cite{jiang2020unified,grapiglia2022adaptive},  universal generalizations for the setting when the $p$-th derivative is H\"older continuous \cite{grapiglia2020tensor,song2021unified,doikov2022super}, versions with inexact minimization of the Taylor model \cite{doikov2020inexact,grapiglia2021inexact}, tensor methods for finding approximate stationary points of convex functions \cite{grapiglia2020tensorNorm,dvurechensky2019near-optimal}.
The ideas described above, especially the Monteiro--Svaiter accelerated proximal point method, turned out to be productive also in the areas not directly related to tensor methods, see non-trivial examples in \cite{bubeck2019complexity,bullins2020highly,carmon2021thinking}. Modern second-order and third-order methods demonstrate also their efficiency in Data Science and Machine Learning applications \cite{daneshmand2021newton,agafonov2021accelerated,dvurechensky2021hyperfast,bullins2021stochastic}. More details of such applications are given in Section~\ref{stoch}.

\section{Notation and problem statement}
The following optimization problem is considered
\begin{equation}
\label{eq1}
\min\limits_{x \in \R^d }\{ F\left( x \right):=f\left( x \right)+g\left( x \right)\} ,
\end{equation}
where $f$ and $g$ are convex functions. Denote $x^*$ -- the solution of the problem \eqref{eq1}. If the solution is not unique, let  $x^*$ be a solution that is the closest to starting point $x_0$ in $2$-norm.

Denote $\|\,\cdot\,\|$ the Euclidean $2$-norm in $\R^d$, $$D^k f(x)[h]^k = \sum_{i_1,...,i_d \ge 0:\,\, \sum_{j=1}^d i_j = k} \frac{\partial^k f(x)}{\partial x_1^{i_1}...\partial x_d^{i_d}}h_1^{i_1} \cdot...\cdot h_d^{i_d},$$  $$\|D^k f(x)\| = \max_{\|h\|\le 1} \left\|D^k f(x)[h]^k\right\|.$$
Assume that $f$ has Lipschitz derivatives of order $p$ ($p \ge 1$):
\begin{equation}
    \|D^p f(x)- D^p f(y)\|\leq L_{p,f}\|x-y\|, \;x,y \in \R^d.
    \label{def_lipshitz}
\end{equation}
Here and below (see e.g. \eqref{unif_conv}) it is considered that $x,y \in \R^d$ belongs to the ball in $2$-norm centred at $x^{*}$ with the radius $O(\|x_0 - x^*\|)$ \cite{nesterov2021inexact}.

The $p$-th Taylor polynomial of $f$ is defined as
\begin{equation}
    \Omega_{p}(f,x;y)=f(x)+\sum_{k=1}^{p}\frac{1}{k!}D^{k}f(x)\left[ y-x \right]^k, \; y\in \R^d.
    \label{eq_taylor}
\end{equation}
Note that from \eqref{def_lipshitz} it follows  \cite{nesterov2021implementable} that
\begin{equation}
   \left|f(y)-\Omega_{p}(f,x;y)\right| \leq \frac{L_{p,f}}{(p+1)!}\|y-x\|^{p+1}.
    \label{eq_sumup}
\end{equation}

$F$ satisfies \textit{$r$-growth condition} ($p+1 \geq r \geq 2$) (see \cite{poljak1982sharp}) with constant $\sigma_r > 0$ iff
\begin{equation}\label{unif_conv}
    F(x) - F(x_*) \ge \sigma_r \|x-x^*\|^r, \; x \in \R^d.\end{equation}

\section{Optimal Tensor methods}\label{optimal}
The following algorithm is taken from \cite{gasnikov2021accelerated}, for $g\equiv 0$ see \cite{bubeck2019near}.
\begin{algorithm} [h!]
\caption{Monteiro--Svaiter--Nesterov \\ MSN($x_0$,$f$,$g$,$p$,$H$,$K$)}\label{alg:highorder}
	\begin{algorithmic}[1]
		\STATE \textbf{Input:} $p \ge 1$, $f : \R^d \rightarrow \R$, 
		$g : \R^d \rightarrow \R$, $H > 0$. 
		\STATE $A_0 = 0, y_0 = x_0$.
		\FOR{ $k = 0$ \TO $k = K- 1$}
		\STATE Find $\lambda_{k+1} > 0$ and $y_{k+1}\in \R^d$ such that
		\[
\frac{1}{2} \leq \lambda_{k+1} \frac{H \|y_{k+1} - \tilde{x}_k\|^{p-1}}{p!}  \leq \frac{p}{p+1} \,, \text{ where}
\]
\begin{equation}
\label{prox_step}
\hspace{-2em} y_{k+1} = \argmin_{y\in \R^d} \left\{
\Omega_{p}(f,\tilde{x}_k;y)+g(y) +\frac{H}{(p+1)!}\|y-\tilde{x}_k\|^{p+1} \right\} \,,
\end{equation}
		\[
a_{k+1} = \frac{\lambda_{k+1}+\sqrt{\lambda_{k+1}^2+4\lambda_{k+1}A_k}}{2} 
\text{ , } 
A_{k+1} = A_k+a_{k+1}
\text{ , } 
\]
\[
\tilde{x}_k = \frac{A_k}{A_{k + 1}}y_k + \frac{a_{k+1}}{A_{k+1}} x_k. 
		\]
		\STATE $x_{k+1} := x_k-a_{k+1} \nabla f(y_{k+1}) - a_{k+1}\nabla g(y_{k+1})$.
		\ENDFOR
		\RETURN $y_{K}$ 
	\end{algorithmic}	
\end{algorithm}

\begin{theorem} \label{theoremCATD} \textup{\cite{gasnikov2021accelerated}}
Let $y_k$ be an output point of Algorithm~\ref{alg:highorder} MSN($x_0$, $f$, $g$, $p$, $H$, $k$) after $k$ iterations, when $p\geq 1$ and $H\ge (p+1)L_{p,f}$. Then
 \begin{equation} \label{speedCATD}
 F(y_k) - F(x^{\ast}) \leq \frac{c_p H R^{p+1}}{k^{\frac{3p +1}{2}}} \,,
 \end{equation}
where $c_p = 2^{p-1} (p+1)^{\frac{3p+1}{2}} / p!$,
$R=\|x_0 - x^{\ast}\|$. 

Moreover, when $p \ge 2$ for $\varepsilon$: $F(y_k) - F(x_{\ast}) \leq \varepsilon$ it is required to solve auxiliary problem \eqref{prox_step}, to find  $(\lambda_{k+1},y_{k+1})$ with proper accuracy, $O\left(\ln\left(\varepsilon^{-1}\right)\right)$ times.
\end{theorem}

Note that Theorem~\ref{theoremCATD} is true also when $H\ge 2L_{p,f}$ (independently of $p \ge 1$). It can be derived from \eqref{eq_sumup}. The condition $H\ge (p+1)L_{p,f}$ was used only because it guarantees the convexity of auxiliary problem \eqref{prox_step}, see \cite{nesterov2021implementable}. Under the conditions $H\ge (p+1)L_{p,f}$, $g \equiv 0$ for $p = 1, 2, 3$ there exist efficient ways to solve auxiliary problem \eqref{prox_step}, see \cite{nesterov2021implementable}. For $p = 1$ there exists an explicit formula for the solution of \eqref{prox_step}. For $p = 2, 3$ the complexity of \eqref{prox_step} is almost the same (up to a logarithmic factor in $\varepsilon$) as a complexity of Newton's method iteration  \cite{nesterov2021implementable}, see also Section~\ref{hyperfast}. 
It is important that there is no need to solve \eqref{prox_step} accurately. It is sufficient to find such $\tilde{y}_{k+1}$ that satisfied
\begin{equation}
\label{inexact1}
    \left\|\nabla \widetilde{\Omega}^k(\tilde{y}_{k+1}) \right\| \le \frac{1}{4p(p+1)}\|\nabla F(\tilde{y}_{k+1})\|,
\end{equation}
where
$$\widetilde{\Omega}^k(y):=
\Omega_{p}(f,\tilde{x}_k;y)+g(y) +\frac{H}{(p+1)!}\|y-\tilde{x}_k\|^{p+1}.$$
Such a practical modification slow down the theoretical convergence rate in a factor $12/5$ in the right hand side of \eqref{speedCATD} \cite{kamzolov2020near,grapiglia2021inexact,nesterov2021inexact}.

In May 2022, D. Kovalev et al. \cite{kovalev2022first_} proposed for $g \equiv 0$ an explicit policy for choosing a pair  $(\lambda_{k+1},y_{k+1})$ in Algorithm~\ref{alg:highorder} and modify the stopping criteria \eqref{inexact1} according to  \cite{ivanova2021adaptive}. It allows to solve \eqref{prox_step} on each iteration at average only two times rather $O\left(\ln\left(\varepsilon^{-1}\right)\right)$ times as it was before. The final complexity bound for $p$-order oracle calls in \cite{kovalev2022first_} matched the lower oracle complexity bound from \cite{nesterov2021implementable} (see also \cite{garg2021near}) obtained for the worst-case function  
$$F_p\left(x_1,...,x_d\right) = |x_1|^{p+1} + |x_2 - x_1|^{p+1} + ... + |x_d - x_{d-1}|^{p+1}.$$
In concurrent and independent paper \cite{carmon2022optimal}, the authors also proposed a way of reducing additional logarithmic factors. Note, that approach of \cite{carmon2022optimal} does not require any a priory knowledge of smoothness parameters (including Holder continuity assumption instead of Lipschitz one).

If additionally $F$ satisfies $r$-growth condition \eqref{unif_conv}
then optimal method can be developed based on restarts procedure \cite{nesterov2006cubic,gasnikov2019optimal} -- see Algorithm~\ref{alg:restarts}.

\begin{algorithm} [h!]
\caption{Restarted MSN($x_0$,$f$,$g$,$p$,$r$,$\sigma_r$,$H$,$K$)}\label{alg:restarts}
	\begin{algorithmic}[1]
		\STATE \textbf{Input:}  $F = f + g : \R^d \rightarrow \R$ satisfies $r$-growth condition with constant $\sigma_r$, MSN($x_0$,$f$,$g$,$p$,$H$,$K$).
		\STATE $z_0=x_0$.
		\FOR{$k = 0 $ \TO $K$}
		\STATE $R_k=R_0\cdot 2^{-k}$, 
		\begin{equation}
		\label{numberofrestarts}    
		N_k=\max \left\{ \left\lceil \left(\frac{r c_p H 2^r}{\sigma_r} R_k^{p+1-r}\right)^{\frac{2}{3p+1}}  \right\rceil, 1\right\}.
		\end{equation}
		\STATE  $z_{k+1} := y_{N_k}$, where $y_{N_k}$ -- output of MSN($z_k$,$f$,$g$,$p$,$H$,$N_k$).
		\ENDFOR
		\RETURN $z_{K}$ 
	\end{algorithmic}	
\end{algorithm}
\begin{theorem} \label{theoremRestartCATD} \textup{\cite{nesterov2006cubic,gasnikov2019optimal}}
Let $z_K$ be an output of Algorithm~\ref{alg:restarts} after $K$ restarts. If $H\ge (p+1)L_{p,f}$, $\sigma_r > 0$, then  for $F(z_K) - F(x^*) \le \varepsilon$ it sufficient to solve \eqref{prox_step}:
\begin{equation}\label{r}
    N = \tilde{O} \left(\left( \frac{H R^{p+1-r}}{\sigma_r} \right)^{\frac{2}{3p+1}}\right)
\end{equation}
times, where $\tilde{O}(\,)$ -- means the same as $O(\,)$ up to a
$\ln^{\alpha}\left({\varepsilon^{-1}}\right)$ factor.
\end{theorem}
Everything that was noted after Theorem~\ref{theoremCATD} will also take place in this case.

In particular, when $g \equiv 0$, $p\ge 2$ and $r = 2$ MSN algorithm in Theorem~\ref{theoremRestartCATD} can be replaced by Kovalev's variant of MSN \cite{kovalev2022first_} to obtain 
\begin{equation}\label{sconv}
O\left(\left(\frac{L_{p,f}R^{p-1}}{\sigma_2}\right)^{\frac{2}{3p+1}} + \log \log \left(\frac{\sigma_2^{\frac{p+1}{p-1}}}{L_{p,f}^{\frac{2}{p-1}}\varepsilon} \right) \right)
\end{equation}
$p$-order oracle complexity bound. This upper bound corresponds to the $\sigma_2$-strongly convex case lower bound from \cite{kornowski2020high} and improves \eqref{r} on a logarithmic factor. The dependence on $\varepsilon$ is $\log \log \varepsilon^{-1}$ as it should be locally for tensor methods \cite{gasnikov2017universal,doikov2022local} (see also \cite{nemirovskij1983problem,nesterov2006cubic,agarwal2018lower,doikov2022super}). But \eqref{sconv} describe two regimes: the first term (complexity independent on $\varepsilon$) describe the complexity to reach the vicinity of quadratic convergence, the second term (complexity is $\log \log \varepsilon^{-1}$) describe Newton's type local convergence.

Note that due to the presence of composite term $g$ the described above MSN algorithm and its variations can be used for splitting oracle complexities \cite{kamzolov2020optimal,kovalev2022optimal_sliding}. Namely, assume that $f$ and $g$ have Lipshitz $p$-th order derivatives with different Lipschitz constants. Based on MSN algorithm one can propose a general framework to accelerate tensor methods by splitting computational complexities. As a result, one can get near-optimal oracle complexity for each function in the sum separately for any $p\geq 1$, including the first-order methods. To be more precise, if the near optimal complexity to minimize $f$ is $N_f(\varepsilon)$ and to minimize $g$ is $N_g(\varepsilon)$, then MSN-based sliding algorithm \cite{kamzolov2020optimal,kovalev2022optimal_sliding} requires no more than $\tilde{O}\left(N_f(\varepsilon)\right)$ oracle calls for $f$ and $\tilde{O}\left(N_g(\varepsilon)\right)$ oracle calls for $g$ to minimze $F = f+g$.

Not also that for $p=1$ and $f \equiv 0$ one can take $H > 0$ in MSN and obtain Monteiro--Svaiter accelerated proximal envelop. In the cycle of recent papers \cite{ivanova2021adaptive,carmon2022recapp,kovalev2022optimal_sliding,kovalev2022first,kovalev2022first_} it was shown that this type of proximal (Catalyst-type \cite{lin2018catalyst}) envelop is logarithmic-free due to the well developed stopping rule criteria for the inner (auxiliary) problem. So it opens up different perspectives for improving existing Catalyst-based algorithms, see e.g. \cite{tominin2021accelerated,tian2022acceleration}. 

\section{Superfast acceleration and the structure of the auxiliary problem}\label{hyperfast}

In this section, the computationally efficient solution of the tensor subproblem for $p=3$ is discussed. With the proposed procedure, it is possible to implement third-order methods without the computation and storing of the third-order derivative. At the beginning of the section, it is shown that this subproblem is convex. Then the Bregman-distance gradient method is used as a subsolver for the tensor step. \eqref{prox_step}. 

For this section, it is assumed that $g = 0$ and $p = 3$.
Then $L_3 = L_{3,f}$, the third-order Taylor's polynomial is
\begin{equation}
\begin{gathered}
\Omega(f,x;y)=\Omega_{3}(f,x;y)= f(x)+\nabla f(x)[y-x] \\+ \tfrac{1}{2} \nabla^2 f(x)\left[ y-x\right]^2
+ \tfrac{1}{6} D^{3}f(x)\left[y-x \right]^3.
\end{gathered}
\end{equation}
The regularized third-order model is
\begin{equation}
\label{super_model}
\begin{gathered}
\tilde{\Omega}(f,x;y)= f(x)+\nabla f(x)[y-x] \\+ \tfrac{1}{2} \nabla^2 f(x)\left[ y-x\right]^2
+ \tfrac{1}{6} D^{3}f(x)\left[y-x \right]^3 + \tfrac{H}{24}\|y-x\|^{4}.
\end{gathered}
\end{equation}
The basic step for every tensor method is formulated as
\begin{equation}
\label{super_problem}
\begin{gathered}
x_{k+1} =  \argmin\limits_{y\in \R^d}  \left\lbrace   \tilde{\Omega}(f,x_k;y)\right\rbrace.
\end{gathered}
\end{equation}
This step is a major part of almost every third-order method. Next, it is described how to effectively solve this subproblem without the computation of the third-order derivative and only using second-order information.  

First, it is clarified that the function \eqref{super_model} is convex for $H\geq 3 L_3$.
Note, that if $H\geq L_3$ the model \eqref{super_model} is possibly non-convex because of third-order derivative. It means that the minimizing non-convex subproblem is harder than minimizing original convex problem. In 2018 Yu.~Nesterov made the breakthrough in \cite{nesterov2021implementable}. It was shown that if $H\geq 3 L_3$ then the model \eqref{super_model} is convex and high-order methods are implementable.

To give some intuition, a sketch of the convexity proof is presented. The following inequality can be derived from the Lipschitz-continuous condition and the convexity of function $f$.
$$0\preceq \nabla^2 f(y)\preceq \nabla^2\Omega(f,x;y) + \frac{L_{3}}{2}\|y-x\|^{2}I.$$
Now, to finish the proof, one need to choose $H$ such that $$\nabla^2 \tilde{\Omega}(f,x;y) \succeq \nabla^2\Omega(f,x;y) + \frac{L_{3}}{2}\|y-x\|^{2}I. $$
This lead to the crucial detail, that was misleading before
\begin{gather*}
    \nabla^2\lb \tfrac{1}{4} \|x\|^4\rb =  2 x x^T + \|x\|^2 I
    \succeq \|x\|^2 I.
\end{gather*}
So, from the fact that $x$ is a vector and $x x^T$ is a singular matrix, the factor $3$ in the last inequality is losing. Finally, if $H\geq 3L_3$ then the function \eqref{super_model} is convex. For more details one can check Theorem 1 in \cite{nesterov2021implementable}.

Thus, one can move to the effective subsolver of the problem \eqref{super_problem} by Bregman-distance gradient method for relatively smooth functions from \cite{lu2018relatively}.
The main idea is to show that the optimized function $\phi(y)$ is $L_{\rho}$-smooth and $\mu_{\rho}$-strongly convex with respect to some convex function $\rho(y)$ or
$$\mu_{\rho}\nabla^2\rho(y)\preceq \nabla^2 \phi(y) \preceq L_{\rho}\nabla^2\rho(y).$$
Then Bregman-distance gradient method make next steps
\begin{equation*}
    y_{t+1} = \argmin\limits_{y \in \R^d} \lb \la \nabla  \phi(y_t), y - y_t  \ra  +
   L_{\rho}  \beta_{\rho}\ls y_t, y \rs \rb,
\end{equation*}
where
$$\beta_{\rho}\ls y_t, y \rs = \rho(y) - \rho(y_t) - \la \nabla \rho(y_t), y - y_t\ra$$
is a Bregman-distance generated by $\rho(y)$.
Bregman-distance gradient method  converges linearly with condition number $\kappa = \tfrac{L_{\rho}}{\mu_{\rho}}$ and convergence rate 
$$N = \tfrac{L_{\rho}}{\mu_{\rho}} \log \ls \tfrac{L_{\rho} \beta_{\rho}\ls y_0, y_{\ast} \rs}{\varepsilon} \rs.  $$

One can show that the model \eqref{super_model} with $H=6 L_3$ can be optimized as $\phi(y)=\tilde{\Omega}(f,x_k;y)$ by gradient method with Bregman-distance generated by
$$
\rho_{x_k}(y)= \tfrac{1}{2} \nabla^2 f(x_k)\left[ y-x_k\right]^2 + \tfrac{L_3}{4}\|y-x_k\|^{4},
$$
$L_{\rho} = 1 + \tfrac{1}{\sqrt{2}}$, $\mu_{\rho} = 1 - \tfrac{1}{\sqrt{2}}$, and $\kappa= \tfrac{1}{\ls1+\sqrt{2} \rs^2} = O(1)$. This means that this method is very fast and converges with a fixed number of iterations. It is worth noting that for each step, the full hessian for $\beta$ is computed, but the full third-order derivative is unnecessary because one only need the derivative-vector-vector product $D^{3}f(x_k)\left[y-x \right]^2$ to compute $\nabla \tilde{\Omega}(f,x_k;y)$. Derivative-vector-vector product can be efficiently and precisely computed by autogradient or approximated by finite difference. 
To summarize, the tensor method with a convergence rate of $1/k^5$ is obtained for third-order methods but only the first and second-order derivatives are computed.
This approach was firstly proposed in \cite{nesterov2021implementable} and then it was improved in \cite{nesterov2021superfast}. Section 5 from \cite{nesterov2021implementable} is good for a general understanding of this approach. For details and precise proofs, one can read \cite{nesterov2021superfast}.

\section{Tensor methods and stochastic distributed setup}\label{stoch}

As well as first-order methods, modifications of tensor methods are proposed for problems of the finite sum type, stochastic and distributed optimization. 

\subsection{Stochastic Optimization}

In this subsection, the problem 
\eqref{eq1}  with $g=0$ under inexact information on the derivatives of  objective $f$ up to order $p\geq 2$ is considered. In particular, it is motivated by stochastic optimization. In this case, objective $f$ has the following form
\begin{equation}\label{eq:online_intro}
     f(x)=\mathbb{E}_{\xi \sim \mathcal{D}}[f(x ; \xi)],
\end{equation} 
where the random variable $\xi$ is sampled from a distribution $\mathcal{D}$, and an optimization procedure has access only to stochastic realizations of the function $f( x ; \xi)$ via its higher-order derivatives up to the order $p \ge 2$.

The full Taylor expansion of $f$ \eqref{eq_taylor} requires computing all the derivatives up to the order $p$, which can be expensive to calculate. Following the works \cite{ghadimi2017second, agafonov2020inexact,agafonov2023advancing} some approximations $G_{x, i}$ are used for the derivatives $\nabla^i f\left(x\right)$, $i=1,\ldots,p$ to construct an inexact $p$-th order Taylor expansion of the objective: 
\begin{equation}\label{eq:approx_taylor}
    \Phi_{p}(f, x; y)=f\left(x\right)+\sum \limits_{i = 1}^p  \frac{1}{i!}G_{x, i}[y-x]^i. 
\end{equation}
Let us consider the case when approximations $G_{x, i} = G_i (x, \xi)$ are stochastic, namely, they satisfy the following assumption.
\begin{assumption}
    [Unbiased stochastic gradient with bounded variance and bounded variance stochastic high-order derivatives] 
    \label{as:p_ord_stoch}
        For any $x \in \R^d$ stochastic gradient $G_1(x, \xi)$ and stochastic high-order derivatives $G_i(x, \xi), ~ i=2,\ldots, p$ satisfy
        \begin{gather}
                \mathbb{E}[G_1(x, \xi) \mid x] = \nabla f(x), \quad
                \mathbb{E}\left[\|G_1(x, \xi)-  \nabla f(x)\|^2 \mid x\right] \leq \sigma_1^2, \label{eq:p_ord_stoch_grad_def} \\
                \mathbb{E}\left[\|G_i(x, \xi) -  \nabla^i f(x)\|^2 \mid x\right] \leq \sigma_i^2, ~ i=2, \ldots, p \nonumber.        
        \end{gather}
\end{assumption}

Next, stochastic tensor methods under Assumption \ref{as:p_ord_stoch} is described. 

Recall, that exact tensor methods are based on minimisation of tensor model $\tilde{\Omega}(f, x; y)$ (see eq. \eqref{super_model} for $p=3$). For inexact tensor methods the model is constructed in the following way \cite{agafonov2023advancing, agafonov2020inexact}:
\begin{equation*}
\omega^{H, \bar{\delta}}_{p}(f, x; y) \eqdef \phi_{x, p}(y) + \tfrac{\bar{\delta}}{2}\|x-y\|^2 + \textstyle{\sum \limits_{i=3}^p} \tfrac{\eta_i \sigma_i}{i!}\|x-y\|^i + \tfrac{pH}{(p+1)!}\|x-y\|^{p+1},
\end{equation*}
where $\bar{\delta} > 0$ is regularization parameter dependent on $\sigma_1, \sigma_2$,  and $\eta_i$ are some constants: $\eta_i > 0, ~ 3 \leq i \leq p$. 

This model satisfies two main conditions (\cite{agafonov2020inexact}, Theorem 1): 
\begin{itemize}
    \item Model $\omega^{H, \bar{\delta}}_{p}(f, x; y)$ is a global upper bound for the function $f$:
    \begin{equation*}\label{eq:model_major_simple}
            f(y) \leq \omega^{H, \bar{\delta}}_{p}(f, x; y), \; x,y \in \mathbb{R}^n.
        \end{equation*} 
    \item Model $\omega^{H, \bar{\delta}}_{p}(f, x; y)$ is convex.
\end{itemize} 

At each step of the tensor algorithms, one need to find\\
$u \in \argmin_{y \in \R^d} \omega^{H, \bar{\delta}}_{p}(f, x; y)$. However, finding the respective minimizer is a separate challenge. Instead of computing the exact minimum, one can find a point $s \in \R^d$ with a small norm of the gradient.
\begin{definition} \label{def:inexact_sub}
    Denote by $s^{H, \bar \delta, \tau} (x)$ a $\tau$-inexact solution of subproblem, i.e. a point $s \eqdef s^{H, \bar \delta, \tau} (x)$ such that
    \begin{equation*}
        \|\nabla \omega^{H, \bar{\delta}}_{x}(s)\| \leq \tau.
    \end{equation*}
\end{definition}

Everything is now prepared to introduce the Accelerated Stochastic Tensor Method (\cite{agafonov2023advancing}, Algorithm 2) listed as Algorithm \ref{alg:inexact_acc_detailed_p_ord}.

\begin{algorithm} [h!]
  \caption{Accelerated Stochastic Tensor Method \cite{agafonov2023advancing}}\label{alg:inexact_acc_detailed_p_ord}
  \begin{algorithmic}[1]
      \STATE \textbf{Input:} $y_0 = x_0$ is starting point; constants $H \geq \frac{2}{p}L_p$; $\eta_i \geq 4$, ~ $3 \leq i \leq p$; 
     starting inexactness $\bar{\delta}_0  \geq 0$;  nonnegative nondecreasing sequences $\{\bar{\kappa}^t_i\}_{t \geq 0}$ for $i = 2, \ldots, p+1$, and
      \vspace{-8pt}
      \begin{equation}\label{eq:p_ord_alphas}
          \alpha_t = \tfrac{p+1}{t + p + 1}, ~~~ A_t = \textstyle{\prod \limits_{j=1}^t(1 -\alpha_j)}, ~~~ A_0 = 1.
      \end{equation}
      \vspace{-10pt}
      \begin{equation*}
        \psi_{0}(x):= \tfrac{\bar{\kappa}_2^{0}+\lambda_0}{2}\|x - x_0\|^2+ \textstyle{\sum \limits_{i=3}^p} \tfrac{\bar{\kappa}_i^{0}}{i!}\|x - x_0\|^{i} .
    \end{equation*}
    \FOR{$t \geq 0$} 
        \STATE 
                \[v_t = (1 - \alpha_t)x_t + \alpha_t y_t, \quad x_{t+1} = S_p^{H, \bar{\delta}_t, \tau}(v_{t})\]
        \STATE Compute 
            \begin{equation*}
            \begin{aligned}
                &y_{t+1}=\arg \min _{x \in \mathbb{R}^{n}}\left\{\psi_{t+1}(x):= \psi_{t}(x)+ \tfrac{\lambda_{t+1} - \lambda_{t}}{2}\|x - x_0\|^2 \right.\\
                &\left.+ \textstyle{\sum \limits_{i = 2}^{p} }
                \tfrac{\bar{\kappa}^{t+1}_i - \bar{\kappa}^{t}_i}{i!}\|x - x_0\|^i
                +\tfrac{\alpha_{t}}{A_{t}} 
                (f(x_{t+1})+\la G_1(x, \xi), x-x_{t+1}\ra) \right\},
                \end{aligned}
            \end{equation*}
    \ENDFOR
  \end{algorithmic}
\end{algorithm}

\begin{theorem}\label{thm:acc_convergence_p_ord}
\textup{\cite{agafonov2023advancing}}~Let $f$ has a Lipschitz derivatives of order $p$ and Assumption \ref{as:p_ord_stoch} hold.
     After $T \geq 1$ iterations of Algorithm \ref{alg:inexact_acc_detailed_p_ord} with parameters
        \begin{equation*}
        \begin{gathered}
         \label{eq:p_ord_params}
            \bk_2^{t} =  \tfrac{2p\alpha_t^2}{A_t}\bar{\delta}_t;  ~
            \bk_{p+1}^{t} = \tfrac{(p+1)^{p+1}}{2}\tfrac{\alpha_t^{p+1}}{A_t}H; \\
            \bk_{i}^{t} = \tfrac{p^{i-1}}{2}\tfrac{\alpha_t^i}{A_t}\eta_i \sigma_i,~ 3 \leq i \leq p;\\
            \lambda_t = \tfrac{\sigma}{R}(t+p+1)^{p+1/2};~
            \delta_t = \sigma_2 + \tfrac{\tau+\sigma_1}{R}(t+p+1)^{3/2}.
        \end{gathered}
        \end{equation*}
         the objective residual can be bounded as follows
        \begin{equation*}
            \begin{aligned}
               \mathbb{E} \left[ f(x_{T}) - f(x^{\ast}) \right]  
               &\leq 
              \tfrac{3(p+1)^3(\tau + \sigma_1)  R}{(T+p+1)^{1/2}} + \textstyle{\sum \limits_{i=3}^p} \tfrac{2(p+1)^{2i-1}}{i!}\tfrac{\sigma_i R^i}{ (T+p+1)^i} + \tfrac{(p+1)^{2(p+1)}}{(p+1)!} \tfrac{HR^{p+1}}{(T+p+1)^{p+1}}.
            \end{aligned}
            \end{equation*}
    \end{theorem}
The Theorem \ref{thm:acc_convergence_p_ord} yields several intriguing outcomes and results.

 First of all, in the Assumption \ref{as:p_ord_stoch}  the unbiasedness of the gradients, but not of the high-order derivatives was assumed. The general case of inexact (and possibly biased) gradients is considered in \cite{agafonov2020inexact}. In this scenario, a slightly different   version of Algorithm \ref{alg:inexact_acc_detailed_p_ord} converges with the rate $O\ls \delta_1R + \textstyle{\sum \limits_{i=2}^p}\frac{\delta_iR^i}{T^i} + \frac{L_pR^{p+1}}{T^{p+1}} \rs$, where $\delta_i$ is the inexactness in $i$-th derivative ($\|(G_{x,i} - \nabla^i f({x}))[y - x]^{i-1}\| \le \delta_i \|y - x\|^{i-1}, i = 1, \ldots, p.$).

Moreover, Algorithm \ref{alg:inexact_acc_detailed_p_ord} also allows for restarts. The total number of iterations of Restarted Accelerated Stochastic Tensor Method (\cite{agafonov2023advancing}, Algorithm~3) to reach desired accuracy $\e:~f(x_t) - f(x^*)\leq \e$ in expectation is
\begin{equation}\label{eq:ACRNM_complexity_}
    O\left( 
    \tfrac{(\tau + \sigma_1)^2}{\mu \e} 
    +
    \ls\sqrt{\tfrac{\sigma_2}{\mu}} + 1\rs\log\tfrac{f(z_0) - f(x^*)}{\e}
    + 
    \textstyle{\sum} _{i=3}^p \left(\tfrac{\sigma_i R^{i-2}}{\mu}\right)^{\frac{1}{i}}
    +  
    \left(\tfrac{L_p R^{p-1}}{\mu}\right)^\frac{1}{p+1} 
    \right).
\end{equation}

One can use the mini-batch Restarted Accelerated Stochastic Tensor Method  to solve problem \eqref{eq:online_intro}. For simplicity let $p=2$ and assume that the subproblem can be solved exactly, i.e. $\tau = 0$. In this approach, the mini-batched stochastic gradient gradients and Hessians are computed as $\tfrac{1}{r_1}\textstyle \sum_{i=1}^{r_1} \nabla f(x, \xi_i), ~\tfrac{1}{r_2}\textstyle \sum_{i=1}^{r_2} \nabla^2 f(x, \xi_i)$ respectively, where $r_1$ and $r_2$ represent the batch sizes. From the convergence estimate \eqref{eq:ACRNM_complexity_}, one can determine the overall number of stochastic gradient computations  $O\ls \tfrac{\sigma_1^2}{\e \mu^{2/3}}\rs$, which is similar to the accelerated SGD method \cite{ghadimi2013optimal}. Interestingly, the number of stochastic Hessian computations scales linearly with the desired accuracy $\e$, i.e., $O\left(\tfrac{\sigma_2}{\mu^{1/3}}\log\frac{1}{\e}\right)$. 

Finally, note that in the case of $p=2$ the Algorithm \ref{alg:inexact_acc_detailed_p_ord} matches the lower bound (\cite{agafonov2023advancing}, Theorem~4.2) ~$\Omega \ls \tfrac{\sigma_1R}{\sqrt{T}}+ \frac{\sigma_2 R^2}{T^2} + \frac{L_2R^3}{T^{7/2}}\rs$ for stochastic second-order methods in gradient and Hessian error terms.

\subsection{Distributed optimization}

In this subsection,  distributed  empirical risk  minimization problem is considered.  In this setup, multiple agents/devices/workers collaborate to solve a machine learning problem by communicating over a central node (server). The goal is to minimize the following finite-sum objective:
\begin{equation}\label{eq:problem}
    \min \limits_{x \in \R^d} \lb f(x):= \tfrac{1}{m}\sum \limits_{j=1}^m f_j(x)\rb,
\end{equation}
 where $f_j(x)$ is the loss function (parametrized by $x \in \R^d$) associated with the data stored on the $j$-th machine and the $\mu_f$-strongly convex function $f(x)$ is the global empirical loss function. Each node has access only to its local data and can communicate with the central server. 
 
Problems of the form \eqref{eq:problem} arise in  distributed machine learning and other applications.
Beyond first-order methods, second-order methods have been developed for such Federated Learning problems. One of the keys to efficient second-order methods for this setting is exploiting the statistical similarity of the data \cite{Zhang2018,dvurechensky2021hyperfast,daneshmand2021newton,bullins2021stochastic, agafonov2021accelerated}. Using  statistical arguments, it is possible to show that if $f_1(x)=\frac{1}{n}\sum_{l=1}^n \ell(x,y_l)$, where $\ell(x,y_l)$ is the individual loss function of each example $l$ stored at the first device chosen to be the central device, then $\|\nabla^2 f(x) - \nabla^2f_1(x)\|_2 \leq \sigma$, where  $\sigma$ is proportional to $\frac{1}{\sqrt{n}}$. Then, defining 
\begin{equation}
    \label{eq:phi_def}
    \rho(x) = \frac{1}{n}\sum_{l=1}^{n} \ell(x; y_l) + \frac{\sigma}{2}\|x\|_2^2,
\end{equation}
 the global objective $f$ becomes strongly convex and smooth relative to $\rho$, i.e., (cf. Section \ref{hyperfast})
$$
\mu_{\rho}\nabla^2\phi(x)\preceq \nabla^2 f(x) \preceq L_{\rho}\nabla^2\phi(x),
$$
where $L_{\rho}=1$, $\mu_{\rho}={\mu_f}/({\mu_f+2\sigma})$, and $\kappa_{\rho} = {L_{\rho}}/{\mu_{\rho}}=1+2\sigma/\mu_f$. Thus, as in Section \ref{hyperfast} one can apply the Bregman-distance gradient method and its <<accelerated>> variant \cite{hendrikx2020statistically}. 

When implemented, these methods require to minimize at master node (the first node)
$$\min_{y}\left\{\la \nabla  f(y_m), y - y_m  \ra  + L_{\phi} \beta_{\rho}\ls y_m, y \rs\right\}.$$ 
The first term is available due to communications and the last term is a sum-type function with $n$ terms according to \eqref{eq:phi_def}. If $n$ is large enough and $d$ is moderate then Hessian calculation time can dominate Hessian inversion time. That allows to use efficiently Tensor methods. In particular Superfast second-order methods mentioned in the Section~\ref{hyperfast}.

This idea was exploited and analysed in \cite{dvurechensky2021hyperfast} which resulted in an efficient distributed second-order solver with good practical performance and in a certain regime total arithmetic operations complexity being better than that of the existing variance reduced first-order algorithms for problem \eqref{eq:problem}. 

The work of A. Agafonov was supported by the Ministry of Science and Higher Education of the Russian Federation (Goszadaniye), No. 075-00337-20-03, project No. 0714-2020-0005.
The work of  A. Gasnikov was supported by the strategic academic leadership program <<Priority 2030>> (Agreement  075-02-2021-1316 30.09.2021).

\section{Conclusion}
In each iteration higher-order (also called tensor) methods minimize a regularized Taylor expansion of the objective function, which leads to faster convergence rates if the corresponding higher-order derivative is Lipschitz-continuous. Recently a series of lower iteration complexity bounds for such methods were proved, and a gap between upper an lower complexity bounds was revealed. Moreover, it was shown that such methods can be implementable since the appropriately regularized Taylor expansion of a convex function is also convex and, thus, can be minimized in polynomial time. Only very recently an algorithm with optimal convergence rate $1/k^{((3p+1)/2)}$ was proposed for minimizing convex functions with Lipschitz $p$-th derivative. For convex functions with Lipschitz third derivative, these developments allowed to propose a second-order method with convergence rate $1/k^5$, which is faster than the rate $1/k^{3.5}$ of existing second-order methods.

\section{MC codes}
\textbf{90C25,65K10,90C15}

\section{Cross‐References}
\begin{itemize}
    \item Automatic Differentiation: Calculation of Newton Steps; Dixon, L.
    \item Computational Optimal Transport ; Dvurechensky P., Dvinskikh D., Tupitsa N., Gasnikov A.
    \item Stochastic gradient descent; Lan G.
    \item Unified analysis of SGD-type methods; Gorbunov E.
    \item Decentralized Convex Optimization over Time-Varying Graphs; Rogozin A., Gasnikov A., Beznosikov A., Kovalev D.
    \item Stochastic Quasi-Newton Scheme; Jalilzadeh, A.
\end{itemize}

\bibliographystyle{plain}
\bibliography{references}

\end{document}